\documentclass[11pt]{amsart}
\usepackage{amsfonts,amsmath,amsthm,amscd,amssymb,latexsym,cite,verbatim,texdraw,floatflt,
caption2}
\RequirePackage{hyperref}
\usepackage[a4paper]{geometry}

\title{Coarse selectors  of graphs    }
\author{ Igor  Protasov}

\address{I.Protasov: Taras Shevchenko National University of Kyiv, Department of Computer Science and Cybernetics, Academic Glushkov pr. 4d, 03680 Kyiv, Ukraine}
\email{i.v.protasov@gmail.com}

\begin{document}
\begin{abstract} 
We consider a connected graph 
$\Gamma$ as a coarse space and prove that $\Gamma$
 admits a 2-selector if and only if $\Gamma$ is 
 either bounded or coarsely equivalent to $\mathbb{N}$ or $\mathbb{Z}$. We apply this result to 
geodesic metric spaces admitting linear orders compatible with  coarse structures.

\end{abstract}
\maketitle



\section{ Introduction and results}


Selectors and orderings of  topological space have a long history with a lot of remarkable results, see \cite{b2}, \cite{b8}.
The investigations of selectors
and orderings
 of  coarse spaces were initiated in 
 \cite{b3}, \cite{b4}, \cite{b5}.


\vspace{3 mm}

Given a class $\mathcal{K}$, of coarse spaces and $X\in \mathcal{K}$, how can one detect whether $X$ admits a linear order, compatible with the coarse structure of $X$? With usage of selector, this question is answered for discrete coarse spaces \cite{b3}, \cite{b4}, finitary coarse spaces of groups and locally finite graphs \cite{b5}. In contrast to  locally finite graphs, an arbitrary unbounded graph needs not to  have a ray, and this is a technical obstacle to characterize all graphs admitting a 2-selector and a linear order,  compatible with coarse structures. 
This paper is to overcome this obstacle. To this end, we develop and apply a new inductive construction of coarse rays and lines.
\vspace{5 mm}

We recall some basic definitions. 
\vspace{5 mm}

Given a set $X$,
a family $\mathcal{E}$  of subsets of $X\times X$ is called a
{\it  coarse structure} on $X$ if
\vspace{5 mm}

\begin{itemize}
\item{} each $E \in \mathcal{E}$  contains the diagonal $\bigtriangleup _{X}:=\{(x,x): x\in X\}$ of $X$;
\vspace{3 mm}

\item{}  if  $E$, $E^{\prime} \in \mathcal{E}$  then  $E \circ E^{\prime} \in \mathcal{E}$  and
$ E^{-1} \in \mathcal{E}$,    where  $E \circ E^{\prime} = \{  (x,y): \exists z\;\; ((x,z) \in E,  \ (z, y)\in E^{\prime})\}$,    $ E^{-1} = \{ (y,x):  (x,y) \in E \}$;
\vspace{3 mm}

\item{} if $E \in \mathcal{E}$ and  $\bigtriangleup_{X}\subseteq E^{\prime}\subseteq E$  then  $E^{\prime} \in \mathcal{E}$.
\end{itemize}
\vspace{5 mm}

Elements $E\in\mathcal E$ of the coarse structure are called {\em entourages} on $X$.

For $x\in X$  and $E\in \mathcal{E}$ the set $E[x]:= \{ y \in X: (x,y)\in\mathcal{E}\}$ is called the {\it ball of radius  $E$  centered at $x$}.
Since $E=\bigcup_{x\in X}( \{x\}\times E[x]) $, the entourage $E$ is uniquely determined by  the family of balls $\{ E[x]: x\in X\}$.
A subfamily ${\mathcal E} ^\prime \subseteq\mathcal E$ is called a {\em base} of the coarse structure $\mathcal E$ if each set $E\in\mathcal E$ is contained in some $E^\prime \in\mathcal E^\prime$.

The pair $(X, \mathcal{E})$  is called a {\it coarse space}  \cite{b9} or  a {\em ballean} \cite{b6}, \cite{b7}.

A coarse space $(X, \mathcal{E})$  is called
 {\it connected} if,  
 for any $x, y \in X$, there exists $E\in \mathcal{E}$ such that $y\in E[x]$.

A subset  $Y\subseteq  X$  is called {\it bounded} if $Y\subseteq E[x]$ for some $E\in \mathcal{E}$
  and $x\in X$.
If  $(X, \mathcal{E})$  is connected then 
the family $\mathcal{B}_{X}$ of all bounded subsets of $X$  is a bornology on $X$.
We recall that a family $\mathcal{B}$  of subsets of a set $X$ is a {\it bornology}
if $\mathcal{B}$ contains the family $[X] ^{<\omega} $  of all finite subsets of $X$
 and $\mathcal{B}$  is closed   under finite unions and taking subsets. A bornology $\mathcal B$ on a set $X$ is called {\em unbounded} if $X\notin\mathcal B$.
A subfamily  $\mathcal B^{\prime}$ of $\mathcal B$ is called a base for $\mathcal B$ if, for each $B \in \mathcal B$, there exists $B^{\prime} \in \mathcal B^{\prime}$ such that $B\subseteq B^{\prime}$.

Each subset $Y\subseteq X$ defines a {\it subspace}  $(Y, \mathcal{E}|_{Y})$  of $(X, \mathcal{E})$,
 where $\mathcal{E}|_{Y}= \{ E \cap (Y\times Y): E \in \mathcal{E}\}$.
A  subspace $(Y, \mathcal{E}|_{Y})$  is called  {\it large} if there exists $E\in \mathcal{E}$
 such that $X= E[Y]$, where $E[Y]=\bigcup _{y\in Y} E[y]$.

Let $(X, \mathcal{E})$, $(X^{\prime}, \mathcal{E}^{\prime})$
 be  coarse spaces. 
 A mapping $f: X \to X^{\prime}$ is called
  {\it  macro-uniform }  if for every $E\in \mathcal{E}$ there
  exists $E^{\prime}\in \mathcal{E}^{\prime}$  such that $f(E(x))\subseteq  E^{\prime}(f(x))$
    for each $x\in X$.
If $f$ is a bijection such that $f$  and $f ^{-1 }$ are macro-uniform, then   $f  $  is called an {\it asymorphism}.
If  $(X, \mathcal{E})$ and  $(X^{\prime}, \mathcal{E}^{\prime})$  contain large  asymorphic  subspaces, then they are called {\it coarsely equivalent.}

Given a coarse  spaces
$(X, \mathcal{E})$, we denote by 
$exp \ X$ the set of all non-empty subsets of $X$ and endow  $exp \ X$ with the coarse structure $exp \ \mathcal{E}$ with the base 
$\{ exp \ E: E\in \mathcal{E} \}$,
where 
$$(A,B)\in exp \ E 
  \Leftrightarrow A \subseteq E[B], \ \ B\subseteq E[A].$$

Let   $(X, \mathcal{E})$ be coarse  space, 
  $\mathcal{F}$ be 
  a non-empty  subspace  of $exp \ X$. 
  A macro-uniform mapping 
 $f: \mathcal{F} \longrightarrow X$ 
 is called 
 an $\mathcal{F}$-{\it selector} 
 of $(X,\mathcal{E})$ if $f(A)\in A$ for each $A\in \mathcal{F}$. In the case 
 $\mathcal{F}= exp \ X$,
 $\mathcal{F}= \mathcal{B}\setminus \{0\} $,
 $\mathcal{F}= [X]^2$ we get a {\it global selector},
 a {\it bornologous selector}  
 and a {\it 2-selector}  respectively.

Every metric $d$ on a set $X$ defines the coarse structure $\mathcal{E}_d$ on $X$ with the base 
 $\{\{ (x,y) : d(x,y)\leq r \} : r> 0 \}$.
 Given a connected graph $\Gamma$, $\Gamma =\Gamma [V]
 $, we denote by $d$ the path metric on the set 
 $V$ of vertices of $\Gamma$ and consider  $\Gamma$ as the 
 coarse space $(V, \mathcal{E}_d)$. 
  We note that two graphs are coarsely equivalent if and only if they are quasi-isometric [1, Chapter 4].

  For locally finite graphs, the following theorem was proved in [5, Theorem 2]. 
  By $\mathbb{N}$ and $\mathbb{Z}$, we denote the graphs on the sets of natural and integer numbers in which two vertices $a, b$ are incident if and only if $|a-b|=1$.

\vspace{7 mm}

 {\bf Theorem 1. } 
{\it For a graph  $\Gamma$,
 the following statements are equivalent
 \vspace{5 mm}
 
 $(i)$  $\Gamma$ admits a bornologous  selector; 
  
 \vspace{5 mm}
 
 $(ii)$  $\Gamma$ admits a 2-selector; 
  
\vspace{5 mm}
 
 $(iii)$  $\Gamma$ is  either 
bounded or coarsely equivalent to 
  $\mathbb{N}$ and $\mathbb{Z}$.}
 
\vspace{7 mm}
  
  Let $(X, \mathcal{E})$ be a coarse space.
We say that a linear order $\leq$ on $X$ is 
{\it compatible with the coarse structure} 
$\mathcal{E}$
if, for every $E\in \mathcal{E}$, there exists $F\in \mathcal{E}$ such that  
$E\subseteq F$ and if $\{x,y\}\in [X]^2$, $x< y$ 
$(y<x)$
 and $y\in X\setminus F[x]$  then $x^\prime <y$ 
$(y<x^\prime)$ 
  for each $x^\prime \in E[x]$.

\vspace{5 mm}

Let $(X, \mathcal{E})$ be a coarse space, $\leq$ be a linear order on $X$. 
 We say that an entourage $E\in \mathcal{E}$ is interval (with respect to $\leq$) if, for each 
 $x\in X$, there exist $a_x , b_x \in X$ such that 
 $a_x \leq x\leq b_x$ and $E[x]= [a_x, b_x]$.
We say that   
$ \mathcal{E}$  is an {\it interval} coarse structure if 
there is a base of $\mathcal{E}$
consisting of interval entourages. Clearly, if 
 $\mathcal{E}$ is interval then $\leq$ is compatible with $\mathcal{E}$.

\vspace{7 mm}

 {\bf Theorem 2. } 
{\it For an unbounded geodesic metric  space $X$, the 
 following statements hold
 \vspace{3 mm}
 
 $(i)$  $X$ admits a 2-selector if and only if $X$ is coarsely equivalent to
 $\mathbb{N}$ and $\mathbb{Z}$;
  
 \vspace{3 mm}
 
 $(ii)$  if $X$ admits a linear order compatible with the coarse structure of $X$ then $X$ is coarsely equivalent to 
  $\mathbb{N}$ and $\mathbb{Z}$.}

\section{ Proof of Theorem 1 }

The implication $(i)\Rightarrow (ii)$ is evident.
To prove $(iii)\Rightarrow (i)$,  
we note that  $\mathbb{N}$ and $\mathbb{Z}$ admits
bornologous selectors. 
In both cases, we put $f(A)= min \ A$.
By Proposition 5 from  \cite{b4}, $\Gamma$ admits a bornologous selector.

\vspace{5 mm}
We prove $(ii)\Rightarrow(iii)$.
Let $f$ be a 
 2-selector of 
 $\Gamma [V]$.
 We define a binary relation $\prec$ on $V$  as follows: $a\prec b$ iff $a\neq b$ and $f(\{a, b\})=a$.

We use  the Hausdorff metric $d_H $
on the set of all
 finite subsets of $V$ defined by
 \vspace{5 mm}
 
 $d_H (A,B)= max \{ d(a,B)$, $ \ d(b, A): a\in A, b\in B\}$,
$d(a,B)= min \{ d(a,b): b\in B\}$
\vspace{5 mm}

\noindent and  note that the coarse structure on 
$[V]^2$ is defined by $d_H$.

Since $f$ is macro-uniform, there exists the minimal natural number $r$ such that if $A,B\in[V]^2$
and $d_H(A, B)\leq 1$ then $d (f(A), f(B))\leq r$.
We fix and use this $r$.

We recall that a sequence of vertices 
$v_0, \dots , v_m$ is a {\it geodesic path} if $d(v_0, v_m)= m$
and $d(a_i, a_{i+1})= 1$ for  each $i\in \{0, \dots , m-1 \}$.
For $v \in V$ and $m\geq 0$, $B(v,m)$ denotes 
$\{u\in V: d(v, u)\leq m\}$.

\vspace{7 mm}

 {\bf Claim 1.} {\it Let $v,a,b \in V$, $p>0$, 
 $d(a,b)\leq p$ and $a,b \in V\setminus B(v, p+r)$.
 If   $a\prec v$
then $b\prec v$.
}

\vspace{5 mm}

We choose a geodesic path  $u_0,\dots, u_k$ from $a$ to $b$. By the assumption, $\{u_0,\dots, u_k\}\cap B(v,r)=\emptyset$.
Hence,  $u_0\prec v$, $ \ u_1\prec v, \dots ,$ 
$ \ u_k\prec v$, so  $b\prec v$.

\vspace{7 mm}

 {\bf Claim 2.} {\it Let $v \in V$, $p>0$ and 
 $z_0, \dots , z_k, \dots , z_m$ 
 be a sequence in $V$ such that 
 $d(z_i,z_{i+1})\leq p$, $i\in \{0, \dots , m-1 \}$
 and $d(v, \{z_0,\dots, z_{m}\}) = d(v, z_k)$.
 Let $v_0, \dots, v_{t}$ be a geodesic path  such  that 
 $v_0 =v$, $v_t = z_k$.
 Assume that  the following statements are satisfied
 \vspace{5 mm}

$(1)$ $ \  d(z_0, z_i)>p+r$,  $i\in \{k, \dots , m \}$;

\vspace{5 mm}

$(2)$ $ \  d(z_m, z_i)>p+r$,  $i\in \{0, \dots , k \}$;
 
 \vspace{5 mm}

$(3)$ $ \  B(z_0, p+r) \cap  \{v_0, \dots , v_t \}=\emptyset$;
 
  \vspace{5 mm}
  
$(4)$ $ \  B(z_m, p+r) \cap  \{v_0, \dots , v_t \}=\emptyset$;
  \vspace{5 mm}
  
  Then $d(v, z_k)\leq p+r$. 
 
}

We suppose the contrary $t> p+r$ and let $z_0 \prec v$.
Applying Claim 1, we get $z_i \prec v$ for every 
$i\in \{0, \dots , m \}$, in particular, $z_m \prec v$.

By $(3)$, $z_0 \prec z_k$.
By  $(4)$, $z_m \prec z_k$.
By $(1)$ and $(3)$, $z_0 \prec z_m$.
By $(2)$ and $(4)$, $z_m \prec z_0$ and we get a contradiction.

\vspace{5 mm}
  
We use 
Claim 2 in the following form. 

\vspace{5 mm}

{\bf Claim 3.} {\it Let 
 $z_0, \dots , z_m $ 
 be a sequence in $V$ such that 
 $d(z_i,z_{i+1})\leq p$, 
 $P= \{ z_0,\dots, z_{m}\}$,  
  $q= 2(r+p)+1$. 
  Let $v\in V$, $ \  d(v, P)>p+r$ and let 
  $ \  d(v, P)=d(v, z_j) $.
 Assume that 
 $B(z_0 , p+r)\cap \{z_{q+1}, \dots, z_{m} \}=\emptyset$ and 
 $B(z_m , p+r)\cap \{z_{0}, \dots, z_{m-q} \}=\emptyset$.
 
  Then either $j\leq q$ or $j\geq m-q$.}
\vspace{7 mm}

 Now we suppose that $\Gamma[V]$ is unbounded and construct a large subset $S$ of $V$ such that $S$ is asymorphic to $\mathbb{N}$ and $\mathbb{Z}$.
 
 In notations of  Claim 3, we put $p=2r+1$ and choose a geodesic path $P_0$

$$y_{4p}, \dots , y_0, \ b_{4p}, \dots , b_1, \  c, \ 
a_{1}, \dots , a_{4p}, \   x_{0}, \dots , x_{4p}. $$

We take $v\in V$ such that $d(v,c)=n$,  $n=16p+1$ a choose a geodesic path   $v_{0}, \dots , v_{n}, \ $
$v_0 = c$, $v_n =v$.

Let $u$ be a vertex in $P_0$ nearest to $v$.
By Claim 3 with $q=3p$, either 
$u\in \{x_{p}, \dots , x_{4p}\}$
  or $u\in \{y_{4p}, \dots , y_{p}\}$. We
  consider the first option and let $c\prec a_r$.
  
  We show that $B(a_{3p}, r)\cap \{v_0, \dots , v_n \}=\emptyset$.
  Assume the contrary and choose a geodesic path 
  $\{t_0, \dots, t_k \}$ from $u$ to $v$.
  Then
  $$c\prec a_r, \ c\prec u, \ a_{3p} \prec u, \ a_{3p} \prec t_1 , \  \dots , \ a_{3p} \prec t_k, \ a_{3p} \prec v_{n-1}, \ 
  \dots , \ a_{3p} \prec c,$$
  but $c\prec a_{r} $ gives $c\prec  a_{3p}  $,
  a contradiction.
  
  The case $a_r \prec c$ is analogical.
  
  We take $j\in \{1, \dots ,n \}$ such that $v_j \in B(a_{3p}, r)$.
  Then $|3p -j| \leq r$, $d(a_j, v_j)\leq 2r$ and 
  $d(a_j, v_{j+1})\leq p$.
  
  We redenote the sequence $v_{j+1}, \dots , v_{n}, \ $
  by $$a_{j+1}, \dots , a_{8p}, x_0, \dots , a_{8p} $$
  and get the set $P_1$
  
  $$y_{4p}, \dots , y_0, \ b_{4p}, \dots , b_1, \  c, \ 
a_{1}, \dots , a_{8p}, \   x_{0}, \dots , x_{8p}. $$

 We suppose that 
 $B(y_{4p}, 3p) \cap , 
\{ a_{j+1}, \dots , x_{8p}\} \neq \emptyset$
 and let $x$ be a point from this intersection. Then 
 $c\prec x$ and $c\prec y_{4p}$,
 but by the choice of $P_0$ and $c\prec a_r$, we get 
 $ y_{4p}\prec c$.
 
 Then, we can apply above arguments with Claim 3 to get 
 $P_2$ and so on. We repeat this procedure as long as possible (at most $\omega + \omega$ times) and run into two cases:
 
 \vspace{5 mm}
 
 {\it Case 1:} either $\{ a_i : 0<i<\omega\}$ is large or $\{ b_i : 0<i<\omega\}$ is large. 
 We consider the first option and define a mapping 
 $f: \{ a_i : 0<i<\omega\}\longrightarrow \mathbb{N}$
 by $f(a_i)=i$.
 Since $d(a_i , c)=i$ and $d(a_i , a_{i+1})\leq p $, 
 $f$ is an asymorphism.

  \vspace{5 mm}
 
 {\it Case 2:}  $\{ a_i, b_i : 0<i<\omega\}$ is large.
 We define a mapping 
  $f: \{c,  a_i, b_i : 0<i<\omega\}\longrightarrow \mathbb{Z}$
 by $f(c)=0$,
  $f(a_i)=i$,
  $f(b_i)=-i$.
  Since every geodesic path from 
 $\{ a_i: 0<i<\omega\}$  to $\{ b_i: 0<i<\omega\}$ 
 meets $B(c, 3p)$, $f$ is an asymorphism, see Claim 5 in \cite{b5}.

\section{ Proof of Theorem 2 }

We recall that a metric space $(X, d)$ is {\it geodesic} if, for any $x,y\in X$, 
there exists an isometric embedding 
$f: [0, d(x,y)]  \ \rightarrow X$ such that
$f(0)=x$, $f(d(x,y))=y$.
\vspace{5 mm}

$(i)$ In light Proposition 5 from 
\cite{b4}, to Apply Theorem 1, it suffices to show that 
$(X,d)$ is coarsely equivalent to some graph $\Gamma$.

We use the Zorn's lemma, to choose a maximal by inclusion subset $V$ of $X$ such that $B(v,1)\cap B(u,1)=\emptyset$ for all distinct $u,v\in V$.
Since $X$
is geodesic, $V$ is large in $X$. 
We consider a graph  $\Gamma[V]$ with  the set of edges 
$E$ defined as follows: $(u,v)\in E$ if and only if there exists $x\in X$ such that 
$B(x,1)\cap B(u,1)\neq\emptyset$,
$B(x,1)\cap B(v,1)\neq\emptyset$.
Then $V$ is asymorphic to $\Gamma[V]$.
\vspace{5 mm}

$(ii)$  By Proposition 2 from 
\cite{b4}, $X$ admits a 2-selector.
By $(i)$, $X$ coarsely equivalent to $\mathbb{N}$ or
$\mathbb{Z}$.

\end{document}